\documentclass{article}
\usepackage[margin=1in]{geometry} 
\usepackage{amsmath,tikz,subcaption}
\usetikzlibrary{decorations.pathreplacing,decorations.markings}

\newtheorem{thm}{Theorem}[section]

\newtheorem{prop}[thm]{Proposition}

\newtheorem{lem}[thm]{Lemma}
\newtheorem{Def}[thm]{Definition}
\newtheorem{rem}[thm]{Remark}

\newtheorem{lemmaa}{Lemma A}
\newtheorem{lemmab}{Lemma B}
\newtheorem{lemmac}{Lemma C}

\newcommand{\be}{\begin{equation}}
\newcommand{\ee}{\end{equation}}
\newcommand{\ben}{\begin{enumerate}}
\newcommand{\een}{\end{enumerate}}

\newcommand{\e}{{\varepsilon}}

\newcommand{\norm}[1]{\left\lVert#1\right\rVert}
\newcommand{\abs}[1]{\left|#1\right|}
\newenvironment{nproof}[1][Proof.]{\noindent\textit{#1}}{\hfill Q.E.D.}

\allowdisplaybreaks

\title{\Large Finsler manifolds with Positive Weighted Flag Curvature}

\author{ Zhongmin Shen\ and \   Runzhong Zhao}

\date{}

\begin{document}

\maketitle

\begin{abstract}

The flag curvature is a natural Finsler extension of the sectional curvature in Riemannian geometry. 
However, there are many  non-Riemannian quantities  which interact with the flag curvature.  
In this paper, we introduce a notion of weighted flag curvature by modifying the flag curvature using the non-Riemannian quantity, $T$-curvature. 
We show that a forward complete open Finsler manifold with positive weighted flag curvature is necessarily 
diffeomorphic to the Euclidean space, and that a compact Finsler manifold with nonnegative weighted flag curvature and strictly convex boundary is diffeomorphic to a Euclidean ball.

\bigskip
\noindent
 {\bf Keywords:}  \\
{\bf MR(2000) subject classification: } 53C60, 53B40
\end{abstract}

\section{Introduction}

The interrelation between local metric properties and the global topology or global geometry of a manifold has long been one of the central topics in 
Riemannian geometry. Prominent results in this field include the Hadamard-Cartan theorem, Myer's theorem, the Gauss-Bonnet theorem, the sphere theorem, the Bishop-Gromov volume comparison theorem, and many more
(see, e.g. \cite{Che}\cite{Pe}). 

We are particularly interested in the Gromoll-Meyer theorem which asserts that a complete open Riemannian manifold $M$ with positive sectional curvature $K >0$ must be diffeomorphic to $\mathbf{R}^n$.
For Finsler manifolds, the flag curvature is a natural extension of the sectional curvature in Riemannian geometry,
so it is natural to ask whether the Gromoll-Meyer theorem still holds for Finsler manifolds with positive flag curvature. 
The study of this problem will lead to a better understanding of the flag curvature. 
However, the problem is rather sophisticated, since the structure of a Finsler manifold is not only controlled by the flag curvature, but also by some non-Riemannian quantities\cite{Sh1}. 
For example, the $S$-curvature introduced in \cite{Sh3}, was used in a combination with the Ricci curvature to prove a Myer-type theorem and a Bishop-Gromov-type 
volume comparison theorem for Finsler manifolds\cite{Oh}\cite{Sh3}. 
It was also used in Wu's construction of a weighted flag curvature from which he proved a comparison theorem for the Laplacian of distance functions on 
Finsler manifolds\cite{Wu}.  

In this paper, we will use another non-Riemannian quantity, known as the $T$-curvature\cite{Sh2}, to modify the flag curvature, and then establish the Gromoll-Meyer 
theorem for Finsler manifolds. The $T$-curvature was introduced in the study of the variation of lengths of curves on Finsler manifolds. A good understanding of the $T$-curvature will allow us to incorporate many 
useful variation techniques in Riemannian geometry into Finsler geometry. Our proof of the theorem demonstrates the power of such techniques. Another situation where the $T$-curvature 
comes into play is in the study of the geometry of hypersurfaces in Finsler manifolds, where it controls the normal curvature of level surfaces of distance functions\cite{Sh2}.

We define a weighted flag curvature $K^{\alpha}$ by modifying the flag curvature using the $T$-curvature in (\ref{WRiem}). We prove the following 
\begin{thm}\label{GM}
	Let $(M, F)$ be forward complete open Finsler manifold. Assume that $K^{\alpha} >0$ for some $\alpha > 0$, then $M$ is diffeomorphic to $R^n$.
\end{thm} 
The weighted flag curvature condition $K^{\alpha}>0$ reduces to the sectional curvature condition $K >0$ when the Finsler metric is Riemannian. 
Thus Theorem \ref{GM}  generalizes the Gromoll-Meyer theorem.

With the techniques developed in the proof of this Gromoll-Meyer-type theorem, we are also able to establish a rigidity result for compact Finsler manifolds with boundary, which partially extends 
Wu's\cite{WuH} and Sha's\cite{Sha} results on compact $p$-convex Riemannian manifolds.
\begin{thm}\label{GM2}
	Let $(M, \partial M,F)$ be a compact Finsler manifold with strictly convex boundary $\partial M$. Assume that $K^\alpha \ge 0$ for some $\alpha > 0$, then $M$ is diffeomorphic to a 
	closed Euclidean ball.
\end{thm}
\bigskip

This notion of weighted flag curvature might have other applications. 
For example, it is easy to prove a Myer-type theorem with $K^\alpha \ge K_0 > 0$; Controlled weighted flag curvature also leads to a Hessian estimate of some modified distance function.
Thus the weighted flag curvature deserves further study, and we expect that a full understanding of this quantity will lead to vast advancement in global Finsler geometry.

\section{Preliminaries}

A \textit{Minkowski norm} $F$ on a vector space $V$ is a nonnegative function satisfying
\ben
	\item[(a)] $F$ is $C^\infty$ on $V\setminus\{0\}$;
	\item[(b)] $F$ is positively homogeneous of degree $1$, in the sense that $F(tv) = tF(v)$ for all $t\ge 0$ and $v\in V$;
	\item[(c)] $F$ is strongly convex, in the sense that the matrix $g_{ij}(v) = \frac12\left[F^2\right]_{v^iv^j}(v) := \frac12\frac{\partial^2}{\partial v^i\partial v^j}\left[F^2\right](v)$ 
	is positive definite for all $v\ne 0$.
\een
A \textit{Finsler metric} $F$ on a manifold $M$ is a continuous function on the tangent bundle $TM$, which is$C^\infty$ on the slit tangent bundle $TM\setminus 0$, and whose restriction on each tangent space
$T_xM$ is a Minkowski norm. Note that given a nowhere vanishing $C^\infty$ vector field $Y$ on $U\subset M$, the family of inner products on tangent spaces $g_Y : (v,w)\mapsto g_{ij}(Y)v^iw^j$ defines a Riemannian metric on $U$.
\bigskip

For a Finsler metric $F$ on an $n$-dimensional manifold $M$, the length of a piecewise $C^\infty$ curve $c: [a,b]\to M$ is 
\[
	L(c) := \int_a^b F(c'(t))dt
\]
and the "distance" $d(p,q)$ from a point $p\in M$ to another point $q\in M$ is defined as the infimum of the lengths of all piecewise $C^\infty$ curves $c: [a,b]\to M$ with $c(a) = p$ and $c(b) = q$.
We shall note that a Finsler metric is in general not reversible, i.e., $F(v)\ne F(-v)$ for a general vector $v\in TM$. As a consequence, the distance $d$ is in general not symmetric.
The distance is realized by geodesics, which are characterized in local coordinates by 
\[
	\frac{d^2\gamma^i}{dt^2} + 2G^i\left(\gamma(t), \frac{d\gamma}{dt}\right) = 0
\]
where 
\[
	G^i(x,y) = \frac14 g^{il}(y)\left\{\left[F^2\right]_{x^ky^l}(y)y^k - \left[F^2\right]_{x^l}(y)\right\}
\]
are called \textit{geodesic spray coefficients}
with $\left(g^{ij}\right)$ being the inverse matrix of $\left(g_{ij}\right)$. It is easy to show that along a geodesic $\gamma$, $F(\gamma')$ is constant. 
A geodesic $\gamma$ is said to be \textit{normal} if $F(\gamma') = 1$. A normal geodesic $\gamma: [a,\infty)\to M$ is called a \textit{geodesic ray} if $d(\gamma(s), \gamma(t)) = t-s$ for all $a \le s \le t$. 
\begin{Def}
	A Finsler manifold $(M,F)$ is said to be \textit{forward complete} if every geodesic $\gamma: (a,b)\to M$ can be extended to a geodesic $\gamma: (a,\infty)\to M$. 
\end{Def}
\bigskip

Using the geodesic spray coefficients $G^i$, we may define 
\[
	N^i_{j}(y) := \frac{\partial G^i}{\partial y^j}(y)
\]
which is sometimes called the \textit{nonlinear connection} in literatures.
The second-order partial derivative of $G^i$, namely $\Gamma^{k}_{ij} := \frac{\partial^2 G^i}{\partial y^i\partial y^j}$ is known as the coefficients of the \textit{Berwald
connection}. Hence the \textit{covariant derivative} $D$ of the Berwald connection is given in locally coordinates by 
\begin{equation}
	D_yU = \left[dU^i(y) + U^j\Gamma^{i}_{jk}(y)y^k\right]\frac{\partial}{\partial x^i}
\end{equation}
where $y\in T_xM$ and $U$ is a $C^\infty$ vector field on a neighborhood of $x$. It is clear from the coordinate expression of $D$ that we only need $U$ to be defined along some curve whose tangent vector at $x$ is $y$. 
We say that a vector field $U$ along a curve $\gamma$ is \textit{parallel} if $D_{\gamma'}U \equiv 0$. It is well-known that for vector fields $U,V$ along a geodesic $\gamma$, 
\[
	\frac{d}{dt}g_{\gamma'}(U,V) = g_{\gamma'}(D_{\gamma'}U, V) + g_{\gamma'}(U, D_{\gamma'}V)
\]
In particular, a parallel vector field $U$ along a geodesic $\gamma$ satisfies 
\[
	\frac{d}{dt}g_{\gamma'}(\gamma', U) = 0,\quad	\frac{d}{dt}g_{\gamma'}(U, U) = 0
\]

The \textit{Riemann curvature} $R_y = R^i_{\;k}(x,y)\frac{\partial}{\partial x^i}\otimes dx^k : T_xM\to T_xM$ 
is given by 
\begin{equation}
	R^i_{\;k}(x,y) = 2\frac{\partial G^i}{\partial x^k} - \frac{\partial^2 G^i}{\partial x^l\partial y^k}y^l + 2G^l\frac{\partial^2G^i}{\partial y^k\partial y^l} - \frac{\partial{G^i}}{\partial y^l}
	\frac{\partial G^l}{\partial y^k}
\end{equation}
in local coordinates.
In the case when  $F$ is a Riemannian metric, $\left(g_{ij}\right)$ depends only on $x\in M$, and we have that $R^i_{\;k}(x,y) = R^{\;i}_{j\;kl}(x)y^jy^l$ with $R^{\;i}_{j\;kl}$ being the
components of the Riemannian curvature tensor. The Finsler extension of the sectional curvature 
\begin{equation}
	K(y,v) := \frac{g_y(R_y(v),v)}{F(y)^2g_y(v,v) - g_y(y,v)^2}
\end{equation}
is called the \textit{flag curvature} of the flag $(y, \mathop{\rm span}\{y,v\})$ (with pole $y$). Note that $K(y,v)$ depends on $y$ and the linear subspace $\mathop{\rm span}\{y,v\}$ only,
since $R_y$ is a linear transformation. 
\bigskip

We say that a function $f: M\to\mathbf{R}$ is \textit{locally Lipschitz} if it is Lipschitz on every compact subset $K\subset M$.
For a $C^\infty$ function $f: M\to\mathbf{R}$, we define the \textit{Hessian} of $f$ by 
\begin{Def} (\cite{Sh2})
\[
	H^2f(v) := \left.\frac{d^2}{dt^2}f(\gamma(t))\right|_{t=0}
\]
is called the \textit{Hessian} of a function $f$ in the direction $v$,
where $\gamma: (-\e, \e)\to M$ is the geodesic with $\gamma'(0) = v$. 
\end{Def}
A related notion is the convexity of functions. 
\begin{Def}
	A continuous function $f: M\to\mathbf{R}$ is said to be strictly convex 
	if for any compact set $K\subset M$, there are constants $\lambda_K, \delta_K > 0$ such that 
	\[
		f(\gamma(-t)) + f(\gamma(t)) - 2f(\gamma(0)) \ge \delta_Kt^2
	\]
	for all normal geodesics $\gamma$ with $\gamma(0)\in K$ and $t\in [-\lambda_K, \lambda_K]$.
\end{Def}
This notion of convexity is justified by the following easy facts:
\begin{lem}
	Let $f: M\to\mathbf{R}$ be a $C^\infty$ function, then $f$ is strictly convex if
	\[
		H^2f(v) > 0
	\]
	for all $v\in TM\setminus 0$.
\end{lem}
\begin{lem}
	Let $I\subset \mathbf{R}$ be an interval and $\gamma: I\to M$ be a normal geodesic. 
	Suppose that $f: M\to\mathbf{R}$ is continuous and strictly convex, then the function $f\circ\gamma : I\to\mathbf{R}$ is 
	strictly convex.
\end{lem}

For an interval $I\subset \mathbf{R}$, $t$ in the interior of $I$, 
and functions $f, \underline{f}: I\to \mathbf{R}$, we say that $\underline{f}$ \textit{supports} $f$ at $t$ if $\underline{f}(t) = f(t)$ and $\underline{f} \le f$ in a neighborhood 
of $t$. Now let $x\in M$, $v\in T_xM$ with $F(v) = 1$ and $\gamma_v: (-\varepsilon, \varepsilon)\to M$ be the normal geodesic with $\gamma_v(0) = x$ and $\gamma_v'(0) = v$.
\begin{lem}\label{SUPPORT}
	Let $M$ be a forward complete Finsler manifold and $f: M\to\mathbf{R}$ be a continuous function. 
		Suppose that for any compact set $K\subset M$, $x\in K$ and $v\in T_xM$ such that $F(v) = 1$, there is a $C^2$ function 
		$\underline{f}: (-\varepsilon,\varepsilon)\to \mathbf{R}$ supporting $f\circ \gamma_v$ at $0$ and $\underline{f}''(0)\ge m_K$, where $m_K > 0$ is a constant depending only on $K$, then 
		$f$ is strictly convex. 
\end{lem}
\begin{nproof}
	Fix some small $\delta > 0$ and take $L := \{x\in M\mid d(K,x) \le \delta\}$, then by the Hopf-Rinow theorem $L$ is a compact set. So there is $\lambda \ge 1$ so that $d(p, q) \le \lambda d(q,p)$ for all 
	$p,q\in L$, and we claim that $L$ contains the set $\{x\in M\mid d(x,K) \le \delta/\lambda\}$. 

	Indeed, suppose that $p\in M\setminus L$ and $q$ is a point of $K$ so that $d(p,q) = d(p, K) =: l$, let $\gamma:[0,l]\to M$ be a minimal normal geodesic connecting $p$ to $q$. The set 
	$S := \{t\in [0,l]\mid 
	\gamma|_{[t,l]} \subset L\}$ is nonempty since $l \in S$. Let $t_0 = \inf S$, so $\gamma|_{(t_0,l]}\subset L$ and by continuity $d(K, \gamma(t_0)) \le \delta$. 
	If $d(K,\gamma(t_0)) < \delta$ then $\gamma(t_0)$ is an interior point of $L$, thus there is $t'< t_0$ so that $\gamma|_{[t',t_0]}$ lies in the interior of $L$, contradicting the choice of $t_0$.
	Thus $d(K,\gamma(t_0)) = \delta$ and $d(p,q) > d(\gamma(t_0), q) \ge d(q,\gamma(t_0))/\lambda \ge d(K, \gamma(t_0))/\lambda = \delta/\lambda$.

	Now choose $\lambda_K > 0$ so that $2\lambda_K < \delta/\lambda$ and for all $x\in K$, $v\in T_xM$ with $F(v) = 1$, the geodesic $\gamma_v: [-2\lambda_K, 2\lambda_K]\to M$ with $\gamma'(0) = v$ exists. 
	It follows that all such geodesics lie in $L$. In particular, for any such geodesic $\gamma_v$ and $t\in \left[-\frac32\lambda_K, \frac32\lambda_K\right]$, 
	there is a $C^2$ function $\underline{f}_t: (t-\varepsilon, t+\varepsilon)
	\to \mathbf{R}$ supporting $f\circ \gamma_v$ at $t$, and $\underline{f}_t''(t) \ge m_L > 0$. Therefore, the function $t\mapsto f\circ\gamma_v(t) - \frac{m_L}{4}t^2$ is convex on 
	$\left[-\frac54\lambda_K, \frac54\lambda_K\right]$, and it follows that for $t\in [-\lambda_K, \lambda_K]$,
	\[
		f(\gamma_v(-t)) + f(\gamma_v(t)) - f(\gamma_v(0)) \ge \frac{m_L}{4}t^2 + \frac{m_L}{4}t^2 = \delta_Kt^2
	\]
	with $\delta_K = \frac{m_L}{2} > 0$.
\end{nproof}
\bigskip

For a Finsler manifold with boundary $(M, \partial M, F)$, let $x\in \partial M$ and $\mathbf{n}\in T_xM$ be the outward pointing unit normal vector so that $g_\mathbf{n}(\mathbf{n}, v) = 0$ for all 
$v\in T_x\partial M$. Consider a vector $v\in T_x\partial M$ and let $\gamma:(-\varepsilon, \varepsilon)\to\partial M$ be the geodesic of $\partial M$ with respect to the induced Finsler metric, 
with $\gamma'(0) = v$. 
We call 
\[
	\Lambda_\mathbf{n}(v) := -g_\mathbf{n}(\mathbf{n}, D_{\gamma'}\gamma'(0))
\]
the \textit{normal curvature} of $\partial M$ in the direction of $v$, where $D$ denotes the covariant derivative with respect to the Finsler metric on $M$.
\begin{Def}
	We say that the boundary $\partial M$ is strictly convex if $\Lambda_\mathbf{n}(v) > 0$ for all $v\in T\partial M\setminus 0$.
\end{Def}

\section{The T-curvature}

For a tangent vector $y\in T_xM\setminus \{0\}$, let $Y$ be a geodesic field, i.e., a $C^\infty$ vector field whose integral curves are geodesics, defined on a neighborhood of $x$ such that $Y_x=y$. Let $\hat{g}:= g_Y$ and 
 $\hat{D}$ denote the Levi-Civita connection of $\hat{g}$.  
\begin{Def}
We call
\begin{equation}
	T_y (v) := g_y ( D_vV - \hat{D}_vV, y)
\end{equation}
the $T$-curvature, 
where $V$ is a vector field with $V_x=v$. 
\end{Def}
In local coordinates, we have
\begin{equation} 
	T_y(v) = y^lg_{kl}(y) \left[ \Gamma^k_{jm}(v) - \Gamma^k_{jm} (y) \right] v^j v^m.
\end{equation}
from which we obtain immediately:
\ben
\item[(a)] $T_{\lambda y} (v) = \lambda T_y (v)$,  $\forall \lambda >0$ and $\forall v\in T_xM\setminus\{0\}$,
\item[(b)] $T_y(\lambda v) = \lambda^2 T_y(v)$,  $\forall \lambda >0$ and $\forall v\in T_xM\setminus\{0\}$,
\item[(c)] $T_y(y)=0$,
\item[(d)] $\displaystyle\lim_{v\to 0}T_y(v) = 0$.
\een

Denote by $\gamma_y$ the geodesic with $\gamma_y'(0)=y$. Let $E(t)$ be a parallel extension of $v$ along $\gamma_y$.
Define $\dot{T}_y(v)$ by 
\[ \dot{T}_y (v) := \left.\frac{d}{dt} \left[ T_{\gamma_y'(t)} (E(t)) \right]\right|_{t=0}.\]

The properties of $\dot{T}$ follow from those of the $T$-curvature.
\ben
\item[(a)] $\dot{T}_{\lambda y} (v) =\lambda^2 \dot{T}_y(v)$, $ \forall \lambda > 0$ and $\forall v\in T_xM\setminus\{0\}$,
\item[(b)] $\dot{T}_y (\lambda v) = \lambda^2 \dot{T}_y(v)$,  $\forall \lambda >0$ and $\forall v\in T_xM\setminus\{0\}$,
\item[(c)] $\dot{T}_y(y) = 0$.
\een 

In practice, $\dot{T}$ will be calculated in local coordinates using a generic extension $V(t)$ of $v$, by
\[ 
	\dot{T}_y (v) = \left.\frac{d}{dt} \left[ T_{\gamma_y'(t)} (V(t)) \right]\right|_{t=0} - 
	\left.\frac{\partial T_{u}(w)}{\partial w^i}\right|_{
		u =\gamma'(0),
		w = v}
		D_{\gamma_y'}V(0)^i.
\]
In particular, if $V(t) = j(t)E(t)$ where $E$ is parallel along a geodesic $\gamma$ and $j$ is a $C^\infty$ function, we have 
\begin{equation}\begin{aligned}\label{TINCOORD}
	\dot{T}_{\gamma'(t)}(V(t)) =& 
	\frac{d}{dt}T_{\gamma'(t)}(V(t)) - 
	\left.\frac{\partial T_{u}(w)}{\partial w^i}\right|_{
		u =\gamma'(t),
		w = V(t)}
	D_{\gamma'}V(t)^i \\
	=& 
	\frac{d}{dt}T_{\gamma'(t)}(V(t)) - 
	\left.\frac{\partial T_{u}(w)}{\partial w^i}\right|_{
		u =\gamma'(t),
		w = V(t)}
	\frac{j'(t)V(t)^i}{j(t)}\\
	=&
	\frac{d}{dt}T_{\gamma'(t)}(V(t)) - 2\frac{j'(t)}{j(t)}T_{\gamma'(t)}(V(t)).
\end{aligned}\end{equation}

\section{The Weighted Flag Curvature}
	Using the $T$-curvature we define a \textit{weighted flag curvature} $K^{\alpha}(y,v)$ by
	\begin{Def} For $y, v\in T_xM\setminus\{0\}$ such that $g_y(y,v) = 0$, set
	\be \label{WRiem}
		 K^{\alpha}(y,v) := \dfrac1{F(y)^2g_y(v,v)}\left[g_y( R_y(v), v)
		 + \dot{T}_y(v) - 
		 \alpha \dfrac{T_y^2 (v)}{g_y(v, v)}\right].
	\ee
	We say that $K^{\alpha} \geq K$ (resp. $K^\alpha > K$) if for all such pair of tangent vectors $y,v$,
	\begin{equation}\label{LBWRCurv}
		K^{\alpha}(y,v) \geq K (\mathrm{resp.\ } > K).
	\end{equation}
	\end{Def}

	\begin{rem}
	\rm
		We shall remark that the curvature $K^\alpha(y,v)$ as defined above depends not only on $y$ and the tangent plane $\mathop{\rm span}\{y,v\}$, but in general also on the (direction of the) 
		vector $v$, in contrast to the flag curvature $K(y,v)$.
	\end{rem}
	\bigskip

	Now we use the following variation formulae derived in \cite{Sh2}
	to give an estimate of the second variation of lengths of curves under the condition $K^\alpha \ge 0$, which will be of repeated use later. 
	\begin{thm}(\cite{Sh2})
		Let $c:[a,b]\to M$ be a normal geodesic, and 
		\[
			H:(-\varepsilon,\varepsilon)\times[a,b]\to M
		\]
		be a $C^\infty$ variation of $c$, with $H(0,t) = c(t)$ and the variation field $V(t):= \left.\frac{\partial H(s,t)}{\partial s}\right|_{s=0}$ along $c$. 
		Denote by $L(s)$ the length of the curve $H(s, \cdot): [a,b]\to M$, we have 
		\[
			L'(0) = \int_a^bg_{c'}(D_{c'}V, c')dt = g_{c'(b)}(V(b), c'(b)) - g_{c'(a)}(V(a), c'(a))
		\]
		\[\begin{aligned}
			L''(0) =& \int_a^b\left[g_{c'}(D_{c'}(V^\bot), D_{c'}(V^\bot))) - g_{c'}(R_{c'}(V^\bot), V^\bot)\right]dt\\
			&+\left[F(V(b))^2g_{c'(b)}(\kappa_b(0),c'(b)) - F^2(V(a))g_{c'(a)}(\kappa_a(0),c'(a))\right]\\
			&+\left[ T_{c'(a)}(V(a)) - T_{c'(b)}(V(b))\right]   
		\end{aligned}\]
	where $V^\bot = V - g_{c'}(c', V)c'$ is the orthogonal component of $V$ relative to the span of $c'$, and $\kappa_t$ is the geodesic curvature of the curve 
	$H(\cdot, t): (-\varepsilon, \varepsilon)\to M$, given by 
	\[
		\kappa_t(s) = \frac{1}{F\left(\frac{\partial H}{\partial s}\right)^2} \left[\frac{\partial^2{H^i}}{\partial s^2} + 2 G^i\left(\frac{\partial H}{\partial s}\right)\right]\frac{\partial}{\partial x^i}
	\]
	in local coordinates. 
	\end{thm}
	\bigskip

	Let $\gamma:(-\varepsilon,\varepsilon)\to M$ be a normal geodesic, and $p\in M$ be a point not lying on $\gamma$. 
	Choose a minimal normal geodesic $c:[0,l]\to M$ connecting $\gamma(0)$ to $p$, let $j,k:[0,l]\to\mathbf{R}$ be $C^\infty$ functions with $j(0) = 1, k(0) = g_{c'(0)}(\gamma'(0), c'(0)), k(1) = 0$, 
	and $E(t)$ be a parallel vector field along $c$
	with $E(0) = \gamma'(0)^\bot = \gamma'(0) - g_{c'(0)}(c'(0),\gamma'(0))c'(0)$, the orthogonal component of $\gamma'(0)$ relative to the span of $c'(0)$, with respect to the inner product $g_{c'(0)}$.
	Let $c_0: (-\varepsilon, \varepsilon)\to M$ be a curve with $c_0(0) = p$ and $c_0'(0) = j(1)E(1)$.

	We consider the variation $H: (-\varepsilon, \varepsilon)\times [0,l]\to M$ of $c$ for which $H(s,0) = \gamma(s), 
		H(s,l) = c_0(s), H(0,t) = c(t)$ and 
	the variation field $V(t)  = j(t)E(t) + k(t)c'(t)$. Note that $V(0) = \gamma'(0)$ and $V(1) = c_0'(0)$ as desired, by our choices of the functions 
	$j$ and $k$.

	\begin{figure}[h]
	\begin{center}	
	\begin{tikzpicture}
	\begin{scope}[decoration={markings,mark=at position 0.5 with{\arrow{latex}}}]
	    \coordinate (O) at (0,0);
	    \coordinate (A) at (2,-1);
	    \coordinate (B) at (-2,1);
	    \coordinate (E) at (1.4,-0.7);
	    \coordinate (F) at (-1.4,0.7);
	    \coordinate (P) at (10,2);
	    \coordinate (C) at (10.2,0.4);
	    \coordinate (D) at (10.7,3.6);
	    \coordinate (G) at (9.96, 1.4);
	    \coordinate (H) at (10.16,2.6);

	    \draw[postaction={decorate}] (O) to[out=20,in=175] (P);
	    \draw[] (A) to (E);
	    \draw[] (E) to (O);
	    \draw[postaction={decorate}] (O) to (F);
	    \draw[] (F) to (B);
	    \draw[] (C) to [in=-90,out=120](G);
	    \draw[] (G) to [in=-95,out=90](P);
	    \draw[postaction={decorate}] (P) to [out=85,in=-110](H);
	    \draw[] (H) to [out=70,in=-130](D);
	    \node[label={[label distance=0cm]0:$\scriptscriptstyle \gamma(0) = c(0)$}] at (O) {}; 
	    \node[label={[label distance=-0.2cm]135:$\scriptscriptstyle c_0(0) = c(l) = p$}] at (P) {}; 
	    
	    \draw[] (9.8,2.01) -- (9.79,1.81) -- (9.99,1.8);
	    
	    \draw[dashed,postaction={decorate}] (F) to[out=18,in=175] (H);
	    \draw[dashed,postaction={decorate}] (E) to[out=25,in=180] (G);

        \node[label={$\scriptstyle c_0$}] at (10.2,1.6) {};
        \node[label={$\scriptstyle \gamma$}] at (-0.4,-0.4) {};
        \node[label={$\scriptscriptstyle c = H(0,\cdot)$}] at (5,1.5) {};
        \node[label={$\scriptscriptstyle H(\varepsilon,\cdot)$}] at (4,2.1) {};
        \node[label={$\scriptscriptstyle H(-\varepsilon,\cdot)$}] at (6,0.4) {};
        
        \draw[fill] (O) circle (0.03);
        \draw[fill] (P) circle (0.03);
	\end{scope}
	\end{tikzpicture}
	\caption{The variation $H$ under consideration. Here $\gamma$ is a normal geodesic, and $c_0$ is a curve with $c_0'(0)$ orthogonal to $c'(l)$ with respect to $g_{c'(l)}$.}
    \end{center}
	\end{figure}

	\begin{lem}\label{SECVAR}
		In the above setting, suppose that $K^\alpha \ge 0$ on $M$, and $K^\alpha\ge K_0\ge 0$ on the compact set $c({[0,l_0]})$ for some $l_0\in [0,l]$, then
		\[
			L'(0) = -g_{c'(0)}(\gamma'(0), c'(0)),
		\]
		\[
			L''(0) \le \left(1+\frac1\alpha\right)\beta^2\int_0^lj'(t)^2dt - K_0\beta^2\int_0^{l_0}j(t)^2dt
			+ F(V(l))^2g_{c'(l)}(\kappa_b(0),c'(l)) + \left[T_{c'(0)}{\gamma'(0)} - T_{c'(0)}(\gamma'(0)^\bot)\right],
		\]
		where $\beta = \sqrt{g_{c'(0)}(\gamma'(0)^\bot,\gamma'(0)^\bot)} = \sqrt{g_{c'}(E,E)}$.
	\end{lem}

	\begin{nproof}
		By assumption $g_{c'(l)}(c_0'(0), c'(l)) = 0$, so $L'(0) = -g_{c'(0)}(\gamma'(0), c'(0))$. On the other hand,
		observe that $V^\bot(t) = j(t)E(t)$, hence $D_{c'}V^\bot(t) = j'(t)E(t)$ and when $\beta \ne 0$,
		\[\begin{aligned}
			L''(0) =& \int_0^l \left[j'(t)^2g_{c'(t)}(E(t),E(t)) - j(t)^2g_{c'(t)}(R_{c'(t)}(E(t)), E(t))\right]dt 
			- \int_0^l \frac{d}{dt} T_{c'(t)}(V^\bot(t))dt \\& + F(V(l))^2g_{c'(l)}(\kappa_l(0),c'(l))+  \left[T_{c'(0)}(\gamma'(0)) - T_{c'(0)}(\gamma'(0)^\bot)\right] \\
			= & \int_0^l \left[j'(t)^2g_{c'(t)}(E(t),E(t)) - j(t)^2g_{c'(t)}(R_{c'(t)}(E(t)), E(t))\right]dt 
			- \int_0^l \left[\dot{T}_{c'(t)}(V(t)^\bot) + 2\frac{j'(t)}{j(t)}T_{c'(t)}(V(t)^\bot)\right]dt 
			  \\& + F(V(l))^2g_{c'(l)}(\kappa_l(0),c'(l))+  \left[T_{c'(0)}(\gamma'(0)) - T_{c'(0)}(\gamma'(0)^\bot)\right] \\
			= & \int_0^l j'(t)^2\beta^2dt - \int_0^l\left[j(t)^2g_{c'(t)}(R_{c'(t)}(E(t)), E(t)) 
			+ j(t)^2\dot{T}_{c'(t)}(E(t)) + 2j'(t)j(t)T_{c'(t)}(E(t))\right]dt 
			  \\& + F(V(l))^2g_{c'(l)}(\kappa_l(0),c'(l))+  \left[T_{c'(0)}(\gamma'(0)) - T_{c'(0)}(\gamma'(0)^\bot)\right] \\
			= & \int_0^l j'(t)^2\beta^2dt - \int_0^l j(t)^2\beta^2K^\alpha(c'(t),E(t))dt - \int_0^l \left[ j(t)^2 \frac{\alpha T_{c'(t)}(E(t))^2}{\beta^2} + 2j'(t)j(t)T_{c'(t)}(E(t))\right]dt 
			  \\& + F(V(l))^2g_{c'(l)}(\kappa_l(0),c'(l))+  \left[T_{c'(0)}(\gamma'(0)) - T_{c'(0)}(\gamma'(0)^\bot)\right] \\
			= & \int_0^l j'(t)^2\beta^2dt - \int_0^l j(t)^2\beta^2K^\alpha(c'(t),E(t))dt - \int_0^l \alpha\left[\frac{j(t)T_{c'(t)}(E(t))}{\beta} + \frac{\beta j'(t)}{\alpha}\right]^2dt 
			+ \int_0^l \frac1\alpha \beta^2j'(t)^2dt
			  \\& + F(V(l))^2g_{c'(l)}(\kappa_l(0),c'(l))+  \left[T_{c'(0)}(\gamma'(0)) - T_{c'(0)}(\gamma'(0)^\bot)\right] \\
			\le& \left(1+\frac1\alpha\right)\beta^2\int_0^lj'(t)^2dt - K_0\beta^2\int_0^{l_0}j(t)^2dt \\
			&+ F(V(l))^2g_{c'(l)}(\kappa_l(0),c'(l)) + \left[T_{c'(0)}{\gamma'(0)} - T_{c'(0)}(\gamma'(0)^\bot)\right],
		\end{aligned}\]
		where we have used (\ref{TINCOORD}) for the second equality.

		The inequality holds trivially when $\beta = 0$ and we are done.
	\end{nproof}
	\bigskip

	The following two cases are of particular interest to us:		
		\begin{enumerate}
			\item[(a)] when $j(t) \equiv 1$ we have 
				\begin{equation}
					L''(0) \le -K_0l_0\beta^2 + F(V(l))^2g_{c'(l)}(\kappa_l(0),c'(l)) + \left[T_{c'(0)}{\gamma'(0)} - T_{c'(0)}(\gamma'(0)^\bot)\right];
				\end{equation}
			\item[(b)] when $j(t) = 1 - \frac{t}{l}$ we have $V(l) = 0$, hence
				\begin{equation}
					L''(0) \le \left(1+\frac1\alpha\right)\frac{\beta^2}{l} - \frac14K_0l_0\beta^2 + \left[T_{c'(0)}{\gamma'(0)} - T_{c'(0)}(\gamma'(0)^\bot)\right].
				\end{equation}
		\end{enumerate}

    \begin{figure}[h]
	\begin{center}
	\begin{subfigure}{0.4\textwidth}
	\begin{center}	
	\begin{tikzpicture}
	    \coordinate (O) at (0,0);
	    \coordinate (A) at (-3,-1);
	    \coordinate (B) at (3,1);
	    \coordinate (P) at (0,5);
        \coordinate (E) at (2.2,4.9);
        \coordinate (F) at (-2.2,5);
        
        \coordinate (M) at (-1.2,-0.4);
        \coordinate (N) at (-1.2,5);
        
        \coordinate (A1) at (0,0.5);
        \coordinate (A2) at (0,1);
        \coordinate (A3) at (0,1.5);
        \coordinate (A4) at (0,2);
        \coordinate (A5) at (0,2.5);
        \coordinate (A6) at (0,3);
        \coordinate (A7) at (0,3.5);
        \coordinate (A8) at (0,4);
        \coordinate (A9) at (0,4.5);
        
        \coordinate (B1) at (-1.2,0.14);
        \coordinate (B2) at (-1.2,0.68);
        \coordinate (B3) at (-1.2,1.22);
        \coordinate (B4) at (-1.2,1.76);
        \coordinate (B5) at (-1.2,2.3);
        \coordinate (B6) at (-1.2,2.84);
        \coordinate (B7) at (-1.2,3.38);
        \coordinate (B8) at (-1.2,3.92);
        \coordinate (B9) at (-1.2,4.46);
	    
	\begin{scope}[decoration={markings,mark=at position 0.4 with{\arrow{latex}}}]
	    \draw[postaction={decorate}] (B) to (A);
	\end{scope}
	\begin{scope}[decoration={markings,mark=at position 0.5 with{\arrow{latex}}}]
		\draw[postaction={decorate}] (O) to (P);
    \end{scope}
	\begin{scope}[decoration={markings,mark=at position 0.8 with{\arrow{latex}}}]
	    \draw[postaction={decorate}] (E) to[in=0,out=170] (P);
	    \draw[] (P) to (F);
	\end{scope}
	    \draw[dashed] (M) to (N);
    \begin{scope}[decoration={markings,mark=at position 1 with{\arrow{latex}}}]
        \draw[postaction={decorate}] (A1) to (B1);
        \draw[postaction={decorate}] (A2) to (B2);
        \draw[postaction={decorate}] (A3) to (B3);
        \draw[postaction={decorate}] (A4) to (B4);
        \draw[postaction={decorate}] (A5) to (B5);
        \draw[postaction={decorate}] (A6) to (B6);
        \draw[postaction={decorate}] (A7) to (B7);
        \draw[postaction={decorate}] (A8) to (B8);
        \draw[postaction={decorate}] (A9) to (B9);
    \end{scope}
        
        \node[label={$\scriptstyle c_0$}] at (0.5,4.9) {};
        \node[label={$\scriptstyle \gamma$}] at (0.5,-0.5) {};
        \node[label={[align=left]$\scriptscriptstyle V(t) $\\$\scriptscriptstyle= E(t) + k(0)\frac{l-t}{l}\gamma'(t)$}] at (-2.4,2.5) {};;
        \node[label={$\scriptstyle c$}] at (0.3,2.3) {};
        
        \draw[fill] (O) circle (0.03);
        \draw[fill] (P) circle (0.03);
        
	\end{tikzpicture}
	\caption{$j(t)\equiv 1, k(t) = k(0)\frac{l-t}{l}$.}
    \end{center}
    \end{subfigure}\hspace{5em}
    \begin{subfigure}{0.4\textwidth}
	\begin{center}	
	\begin{tikzpicture}
	    \coordinate (O) at (0,0);
	    \coordinate (A) at (-3,-1);
	    \coordinate (B) at (3,1);
	    \coordinate (P) at (0,5);
        \coordinate (E) at (2.2,4.9);
        \coordinate (F) at (-2.2,5);
        
        \coordinate (M) at (-1.2,-0.4);
        \coordinate (N) at (-1.2,5);
        
        \coordinate (A1) at (0,0.5);
        \coordinate (A2) at (0,1);
        \coordinate (A3) at (0,1.5);
        \coordinate (A4) at (0,2);
        \coordinate (A5) at (0,2.5);
        \coordinate (A6) at (0,3);
        \coordinate (A7) at (0,3.5);
        \coordinate (A8) at (0,4);
        \coordinate (A9) at (0,4.5);
        
        \coordinate (B1) at (-1.08,0.14);
        \coordinate (B2) at (-0.96,0.68);
        \coordinate (B3) at (-0.84,1.22);
        \coordinate (B4) at (-0.72,1.76);
        \coordinate (B5) at (-0.60,2.3);
        \coordinate (B6) at (-0.48,2.84);
        \coordinate (B7) at (-0.36,3.38);
        \coordinate (B8) at (-0.24,3.92);
        \coordinate (B9) at (-0.12,4.46);
	    
	\begin{scope}[decoration={markings,mark=at position 0.4 with{\arrow{latex}}}]
	    \draw[postaction={decorate}] (B) to (A);
	\end{scope}
	\begin{scope}[decoration={markings,mark=at position 0.5 with{\arrow{latex}}}]
		\draw[postaction={decorate}] (O) to (P);
    \end{scope}
	    \draw[dashed] (M) to (P);
    \begin{scope}[decoration={markings,mark=at position 1 with{\arrow{latex}}}]
        \draw[postaction={decorate}] (A1) to (B1);
        \draw[postaction={decorate}] (A2) to (B2);
        \draw[postaction={decorate}] (A3) to (B3);
        \draw[postaction={decorate}] (A4) to (B4);
        \draw[postaction={decorate}] (A5) to (B5);
        \draw[postaction={decorate}] (A6) to (B6);
        \draw[postaction={decorate}] (A7) to (B7);
        \draw[postaction={decorate}] (A8) to (B8);
        \draw[postaction={decorate}] (A9) to (B9);
    \end{scope}
        
        \node[label={$\scriptstyle c_0(t)\equiv c_0(0) = c(l)$}] at (0.3,4.9) {};
        \node[label={$\scriptstyle \gamma$}] at (0.5,-0.5) {};
        \node[label={[align=left]$\scriptscriptstyle V(t) $\\$\scriptscriptstyle= \frac{l-t}{l}E(t) + k(0)\frac{l-t}{l}\gamma'(t)$}] at (-2,2.5) {};;
        \node[label={$\scriptstyle c$}] at (0.3,2.3) {};
        
        \draw[fill] (O) circle (0.03);
        \draw[fill] (P) circle (0.03);
        
	\end{tikzpicture}
	\caption{$j(t) = \frac{l-t}{l}, k(t) = k(0)\frac{l-t}{l}$.}
    \end{center}
    \end{subfigure}
    
    \caption{Special cases of varations of interest. For the function $k$ only the values at endpoints, $ k(0) = g_{c'(0)}(\gamma'(0), c'(0)), k(1) = 0$, matter for our estimation.}
    \end{center}
	\end{figure}

\section{Proof of Theorem \ref{GM}}

In this section we prove Theorem \ref{GM}, a generalization of the Gromoll-Meyer theorem. The proof is based on constructing a strictly convex function from the Busemann function on the manifold.

\subsection{The Busemann function}

Let $(M, F)$ be a forward complete Finsler manifold. For a point $p\in M$, denote by $S^+(p, t):= \{x\in M\mid d(p, x) = t\}$ the forward geodesic sphere of radius $t$ centered at $p$.
Put
\[
	b_p^t(x) := t- d(x, S^+(p, t)).
\]

\begin{lem}
\ben
\item[(a)]  $b^t_p(x)$ is bounded
\be
- d(x, p) \leq b^t_p(x ) \leq d(p, x).
\ee 
\item[(b)] 
for any $ d(p, x) \leq t_1 \leq t_2$,
\be
 b_p^{t_1} (x) \geq b_p^{t_2}(x).
\ee
\item[(c)]
for any $x_1, x_2 \in M$,
\be
-d(x_1, x_2) \leq b^t_p(x_1) - b_p^t(x_2) \leq d(x_2, x_1).   \label{eqLb}
\ee
\een
\end{lem}
\begin{nproof} (a)  Take $z\in S^+(p, t)$ such that $d(x,z) = d(x, S^+(p, t))$.
By triangle inequalties
\begin{eqnarray*}
b^t_p(x) & = & t - d(x, z)\\
& = &  d(p, z) - d(x, z) \\
& \leq & d(p, x).\\
b^t_p(x) & = & d(p, z) -d(x, z) \\
& \geq & - d(x, p).
\end{eqnarray*}

(b)
Take $z\in S^+(p, t_2)$ such that $d(x, z)=d(x, S^+(p,t_2))$. Let $c: [0, a]\to M$ be a minimal normal geodesic connecting $x$ to $z$. Let
\[ s_0:= d(x, S^+(p, t_2)) -t_2+t_1.\]
We have $  0 \leq s_0 \leq   d(x, S^+(p, t_2))$. Then
\[  d(p, c(s_0)) \geq t_2 -d(c(s_0), S^+(p, t_2)) = t_2 - d(x, S^+(p, t_2))+s_0 = t_1.\]
Thus $ c(s_0)\in  M\setminus B(p, t_1)$. Now that
\[ d(x, S^+(p, t_1)) \leq d(x, c(s_0))= s_0 =  d(x, S^+(p, t_2))-t_2 +t_1,\]
thus $ b_p^{t_1} (x) \geq b_p^{t_2}(x)$.

(c)  Take $z\in S^+(p, t))$ such that $d(x_1, z) = d(x_1, S^+(p, t))$. Then
\begin{eqnarray*}
b_p^t(x_1) - b_p^t(x_2) & = & d(x_2, S^+(p, t)) - d(x_1, S^+(p, t))\\
& \leq & d(x_2, z) - d(x_1, z)\\
& \leq & d(x_2, x_1).
\end{eqnarray*}
Thus
\[ b_p^t(x_2) - b^t_p(x_1) \leq d(x_1, x_2).\]
and (\ref{eqLb}) follows.
\end{nproof}

\bigskip

Therefore $b_p^t$ converges to a function $b_p$, uniformly on compact sets.

\bigskip

\begin{Def}
The function
\be
	b_p(x):= \lim_{t\to \infty} b^t_p(x).
\ee
is called the {\it Busemann function} at $p$. 
\end{Def}
It follows from (\ref{eqLb}) that $b_p$ is locally Lipschitz:
\[
-d(x_1, x_2) \leq b_p(x_1)-b_p(x_2) \leq d(x_2, x_1).
\]

A crucial step in our proof is to show that the Busemann function is proper. To this end we will need the following definition.

\begin{Def}
	Let $c: [0,\infty)\to M$ be a geodesic ray issuing from $p$, i.e. $c(0) = p$. Define $b_c: M\to \mathbf{R}$ by 
	\[
		b_c(x) := \lim_{t\to\infty} (t - d(x, c(t))),
	\]
	and $b_p: M\to\mathbf{R}$ by 
	\begin{equation}
		\tilde{b}_p(x) := \sup_c b_c(x)
	\end{equation}
	where the supremum is taken over all geodesic rays $c$ issuing from $p$.
\end{Def}

Note that $c(t)\in S^+(p,t)$, hence $b_c(x) \le b_p(x)$. It follows that $\tilde{b}_p(x)\le b_p(x)$, and
\[
	b_p^{-1}(-\infty,a] \subset \tilde{b}_p^{-1}(-\infty,a].
\]
Therefore we only need to show that $\tilde{b}_p$ is proper. We will follow the ideas in \cite{ChGr}.

\begin{lem}\label{TOTCONV}
	Let $M$ be a forward complete Finsler manifold with $K^\alpha > 0$ for some $\alpha > 0$. 
	Denote by $B^-(p,r) := \{x\in M\mid d(x,p)< r\}$ the backward geodesic ball of radius $r$ centered at $p$, and let $c: [0,\infty)\to M$ be a geodesic ray. 
	The set $M\setminus \bigcup_{t> 0}B^-(c(t),t)$ is totally convex, in the sense that for any (normal) geodesic $\gamma: [0,l]\to M$ such that the endpoints 
	$\gamma(0), \gamma(l)\in M\setminus \bigcup_{t> 0}B^-(c(t),t)$, the geodesic $\gamma$ lies entirely in $M\setminus \bigcup_{t> 0}B^-(c(t),t)$.
\end{lem}

\begin{nproof}
	Suppose the contrary, so that there exists a normal geodesic $\gamma: [0,l]\to M$ with 
	$\gamma(0), \gamma(l)\in M\setminus \bigcup_{t> 0}B^-(c(t),t)$ but $\gamma(t) \in B^-(c({s_0}),s_0)$ for some $t\in (0,l)$ and ${s_0} > 0$. 
	Since $\gamma(t) \in B^-(c({s_0}),s_0)$, $d(\gamma(t), c({s_0})) < s_0$. By the triangle inequality, for any $s > s_0$ we have $d(\gamma(t), c(s)) \le d(\gamma(t),c(s_0))
	+d(c(s_0),c(s)) < s_0 + (s - s_0) = s$. 
	
	In particular, for any $s > s_0$, the point on the geodesic $\gamma$ which is closest to $c(s)$ is some interior point $\gamma(t_s)$, and 
	$d(\gamma(t_s),c(s)) < s$, while $d(\gamma(0),c(s)) \ge s, d(\gamma(l),c(s))\ge s$.
	
	Put $p = c(s)$. Without loss of generality, assume that $s > l+1$, so that for any $t\in [0,l]$, $d(\gamma(t),p) \ge d(\gamma(0),p) - d(\gamma(0),\gamma(t)) \ge s - l > 1$.  
	For each $t\in[0,l]$ choose a minimal normal geodesic 
	$c_t:[0,d_t]\to M$ connecting $\gamma(t)$ to $p$, where $d_t = d(\gamma(t),p)$.
	Construct a variation of $c_t$ as in lemma \ref{SECVAR} with $c = c_t$, $l_0 = 1$ and $j(t) = 1 - \frac{t}{l}$, then for every $t$ there is a function $s\mapsto -L_t(s)$ supporting $f:
	s\mapsto -d(\gamma(s),p)$ at $t$ and 
	\[
		-L_t''(t)\ge -\frac{\beta_t^2}{d_t}\left(1+\frac{1}{\alpha}\right) + \beta_t^2\frac{K_0}{4} - \left[T_{c_t'(0)}(\gamma'(t)) - T_{c_t'(0)}(\gamma'(t)^\bot)\right]
	\]
	where $\beta_t = \sqrt{g_{c_t'(0)}(\gamma'(t)^\bot,\gamma'(t)^\bot)}$, and $K_0 > 0$ is a lower bound of the weighted flag curvature 
	in $\overline{B}^+(\gamma([0,l]),1) := \{x\in M\mid d(\gamma([0,l]), x) \le 1\}$, where $\gamma([0,l])$ is the image of $\gamma$ in $M$. 
	Since $d_t = d(\gamma(t),p) \ge s - l$, for large enough $s$ we get that 
	\[
		-L_t''(t)\ge \beta_t^2\frac{K_0}{8} - \left[T_{c_t'(0)}(\gamma'(t)) - T_{c_t'(0)}(\gamma'(t)^\bot)\right].
	\]
	
	Let $\tau_t := \abs{g_{c'_t(0)}(\gamma'(t),c'_t(0))} = \abs{L_t'(t)}$, and choose $\lambda \ge 1$ so that 
	\[
		\frac1{\lambda}\sqrt{g_y(v,v)} \le F(v) \le \lambda \sqrt{g_y(v,v)}
	\]
	for any $x$ on the geodesic $\gamma$, $y\in T_xM\setminus\{0\}$ and $v\in T_xM$. Observe that $\tau_t^2 + \beta_t^2 = g_{c_t'(0)}(\gamma'(t),\gamma'(t))$, hence
	\[
		\frac1{\lambda^2}\le \tau_t^2 + \beta_t^2 \le \lambda^2,
	\]
	and 
	\[
		-L_t''(t)\ge \left(\frac1{\lambda^2} - \tau_t^2\right)\frac{K_0}{8} - \left[T_{c_t'(0)}(\gamma'(t)) - T_{c_t'(0)}(\gamma'(t)^\bot)\right].
	\]
	Since the image of $\gamma$ is compact, there exists $\tilde{\tau} > 0$ so that when $\tau_t \le \tilde{\tau}$, we have $\frac{K_0\tau_t^2}{8} < \frac{K_0}{32\lambda^2}$, and 
	$\abs{T_{c_t'(0)}(\gamma'(t)) - T_{c_t'(0)}(\gamma'(t)^\bot)} < \frac{K_0}{32\lambda^2}$. It follows that
	\[
		-L_t''(t)\ge \frac{K_0}{16\lambda^2} > 0.
	\]
	
	On the other hand, observe that 
	\[
		\beta_t = \sqrt{g_{c_t'(0)}(\gamma'(t)^\bot,\gamma'(t)^\bot)} \le \sqrt{g_{c_t'(0)}(\gamma'(t),\gamma'(t))} \le \lambda F(\gamma'(t)) = \lambda,
	\]
	and 
	\[
		F(\gamma'(t)^\bot) \le \lambda \beta_t \le \lambda^2,
	\]
	we have 
	\[
		-L_t''(t)\ge \beta_t^2\frac{K_0}{8} - 2T_0,
	\]
	where $T_0 = \displaystyle\sup_{t\in [0,l], x = \gamma(t), y,v\in T_xM, F(y) = 1, F(v)\le \lambda^2}\abs{T_y(v)}$. Now let $\varepsilon_0 = \frac{K_0}{16\lambda^2}$ and define 
	$\chi: (-\infty, 0]\to \mathbf{R}$ by 
	\[
		\chi(x) := \int_0^x \exp\left(\frac{2T_0+ \varepsilon_0}{\tilde\tau^2}t\right)dt,
	\]
	so $\chi'(x) = \exp\left(\frac{2T_0+ \varepsilon_0}{\tilde{\tau}^2}x\right) > 0$ and $\chi''(x) = \frac{2T_0+ \varepsilon_0}{\tilde{\tau}^2}\chi'(x)$. Since $\chi$ is (strictly) increasing and $-L_t$ supports $f$ at $t$, 
	$\chi\circ (-L_t)$ supports $\chi\circ f$ at $t$. Moreover, since $-L_t''(t) + \frac{2T_0+ \varepsilon_0}{\tilde{\tau}^2}\tau_t^2 \ge \varepsilon_0 > 0$ for $\tau_t > \tilde\tau$, we have
	\[\begin{aligned}
		(\chi\circ (-L_t))''(t) &= \chi'(-L_t(t)) \left[-L_t''(t) + \frac{2T_0+ \varepsilon_0}{\tilde{\tau}^2}L_t'(t)^2\right] \\
		&= \chi'(-L_t(t)) \left[-L_t''(t) + \frac{2T_0+ \varepsilon_0}{\tilde{\tau}^2}\tau_t^2\right] \ge \chi'(-\sup_t d_t)\varepsilon_0 \ge \chi'(l-s)\varepsilon_0 > 0
	\end{aligned}\]
	by our choice of $\tilde{\tau}$. This means that $\chi\circ f$ is a convex function on $[0,l]$ achieving 
	its maximum at the point $t_s\in (0,l)$ where $\gamma(t_s)$ is closest to $p$ on $\gamma$. This is a contradiction and 
	we conclude that the set $M\setminus \bigcup_{t> 0}B^-(c(t),t)$ is totally convex. 
\end{nproof}
\bigskip
	
	\begin{figure}
	\begin{center}
	\begin{tikzpicture}
	\begin{scope}[decoration={markings,mark=at position 0.5 with{\arrow{latex}}}]
	    \coordinate (O) at (0,0);
	    \coordinate (A) at (0.5,-3);
	    \coordinate (B) at (-0.5,3);
	    \coordinate (P) at (7,1);
	    \coordinate (C) at (2.5,0.482);
	    
	    \coordinate (E) at (-1,-2);
	    \coordinate (F) at (-1.3,1);
	    
	    \coordinate (M) at (0.7,-0.2);

	    \draw[postaction={decorate}] (O) to[out=10,in=180] (P);

	    \draw[] (A) to (B);
	    \draw (C) circle (2.545);
	    
	    \draw[dashed,postaction={decorate}] (M) to (C);
	    \draw[] (E) to[in=-69.4,out=20] (M);
	    \draw[] (M) to[out=110.6,in=-15] (F);
	    
	    \node[label={[label distance=-0.2cm]-30:$\scriptstyle \gamma(t_s) \text{, closest to }c(s)$}] at (M) {}; 
	    \node[label={[label distance=-0.2cm]180:$\scriptstyle \gamma(l)$}] at (F) {}; 
        \node[label={[label distance=-0.2cm]180:$\scriptstyle \gamma(0)$}] at (E) {}; 
        \node[label={[label distance=-0.2cm]90:$\scriptstyle p=c(s)$}] at (C) {}; 
        \node[label={[label distance=-0.2cm]150:$\scriptstyle c$}] at (P) {};
        \node[label={$\bigcup_{t>0}B^-(c(t),t)$}] at (5.5,-2.7) {};
        \node[label={$\scriptstyle B^-(c(s),s)$}] at (2.5,2.1) {};
        \node[label={[align=left]$\scriptstyle \text{Proposed geodesic }\gamma$\\$\scriptstyle\text{ running into }\bigcup_{t>0}B^-(c(t),t)$}] at (-2,-1) {};
        
        \draw[fill] (M) circle (0.03);
        \draw[fill] (F) circle (0.03);
        \draw[fill] (C) circle (0.03);
        \draw[fill] (E) circle (0.03);
	\end{scope}
	\end{tikzpicture}
	\caption{Geodesic $\gamma$ under investigation in the proof of lemma \ref{TOTCONV}.}
	\end{center}
	\end{figure}

\bigskip
\begin{prop}
	On a forward complete Finsler manifold $(M,F)$ with $K^\alpha > 0$, the function $\tilde{b}_p$ is proper for any $p\in M$.
\end{prop}
\begin{nproof}
	Let $c:[0,\infty)\to M$ be a geodesic ray issuing from $p$, and define $c_s:[0,\infty)\to M$ by $c_s(t):= c(s+t)$. 
	Observe now that $b_c^{-1}(-\infty,s] = M\setminus \bigcup_{t> 0}B^-(c_s(t),t) =: P_c(s)$
	and	$\tilde{b}_p^{-1}(-\infty,s]$ is a closed subset of $\bigcap_c P_c(s)$ 
	where the intersection is taken over all normal geodesic rays $c$ issuing from $p$, since $\tilde{b}_p(x) \le s$ implies that $b_c(x) \le s$ for every such ray $c$.
	
	Suppose that $\bigcap_c P_c(s)$ is not compact, so there is a sequence of points $p_k$ in $\bigcap_c P_c(s)$ so that
	$\lim_{k\to\infty}d(p,p_k)= \infty$. For each $k$ choose a minimal normal geodesic $c_k$ connecting 
	$p$ to $p_k$. Up to a subsequence we may assume that $c_k$ converges pointwise to a geodesic ray $c$. 
	Since $\bigcap_c P_c(s)$ is closed and totally convex, every $c_k$ lies entirely in $\bigcap_c P_c(s)$ and so does 
	$c$. But this is a contradiction since $b_c(c(t)) = t$ for all $t > 0$. 
	
	Thus we conclude that $\tilde{b}_p$ is proper. 
\end{nproof}
\bigskip

The next lemma allows for a good estimate of the quantity $b_p(\gamma(-t)) + b_p(\gamma(t)) - 2b_p(\gamma(0))$ for a generic normal geodesic $\gamma$ by applying lemma \ref{SECVAR}.

\begin{lem}\label{georay}
Let $(M, F)$ be a forward complete Finsler manifold and $p\in M$. For any point $q\in M$, there is a geodesic ray $c_q: [0, \infty)\to M$ issuing from $q$ such that 
\ben
\item[(a)] for all $t >0$,
\[ b_p^{q, t} (x) := b_p(q) + t- d(x, c_q(t)).\]
supports $b_p(x)$ at $q$, in the sense that
$ b_p^{q, t} (x) \leq b_p(x)$ for all $x\in M$ and $b_p^{q, t}(q) = b_p(q)$.
\item[(b)] for all $t\geq 0$,
\[ b_p(c_q(t)) = b_p(q)+ t.\]
\een
\end{lem}
\begin{nproof}
Take a sequence $t_n \to \infty$ and a sequence of points $x_n$ in $S^+(p, t_n)$ such that $d(q, x_n)=d(q, S^+(p, t_n))$.  
For each $n$, let $\gamma_n : [0, s_n]\to M$ be a minimal normal geodesic connecting $q$ to $x_n$.  Up to a subsequence we may assume that 
$\gamma_n$ converges pointwise to a geodesic ray $c_q: [0, \infty)\to M$.  
For any sufficiently large $t_n$,
\[ d(q, S^+(p, t_n)) = t+ d(\gamma_n(t), S^+(p, t_n)).\]
Observe that 
\begin{eqnarray*}
b_p(x) - b_p^{q, t}(x) & = &  b_p(x)-b_p(q)-t+d(x, c_q(t))\\
& = &  \lim_{n\to \infty}  \Big \{ [t_n-d(x, S^+(p, t_n))]-[t_n - d(q, S^+(p, t_n))] -t + d(x, c_q(t))   \Big \}\\
&\geq & \lim_{n\to \infty}  \Big \{ - d(x, \gamma_n(t)) +d(x, c_q(t))     \Big \}\\
& \geq &  \lim_{n\to \infty} - d(c_q(t), \gamma_n(t)) = 0. 
\end{eqnarray*}
This proves (a).

On the other hand, for any $s >0$ and sufficiently large $t_n$, 
\[ t_n - d(\gamma_n(s), S^+(p, t_n)) =t_n - [ d(q, S^+(p, t_n)) -s ].\]
\[           -d(c_q(s), \gamma_n(s)) \leq d(\gamma_n(s), S^+(p, t_n)) - d(c_q(s), S^+(p, t_n))\leq d(\gamma_n(s), c_q(s)).\]
Letting $n\to \infty$, we obtain
\[  b_p (c_q(s)) = b_p(q) + s \]
which is (b).
\end{nproof}

\subsection{A Smoothing Theorem}

We have the following theorem for locally Lipschitz strictly convex functions. For the proof we refer the reader to the appendix. An alternate proof was given in \cite{Xia}.
\begin{thm}\label{Smooth}
	Let $M$ be a Finsler manifold, and $f: M\to \mathbf{R}$ be a locally Lipschitz strictly convex function on $M$. Given any $\varepsilon > 0$, there is a locally Lipschitz, strictly convex
	$C^\infty$	function $g: M\to\mathbf{R}$ such that $\abs{g - f} < \varepsilon$ on $M$. 
\end{thm}
Since the set of Morse functions is dense in the space of $C^\infty$ functions, the existence of a proper strictly convex 
Morse function will follow given a proper strictly convex $C^\infty$ function.




\subsection{Proof of Theorem \ref{GM}}

In view of the Morse theory, to prove that $M\cong \mathbf{R}^n$, it now suffices to show that there exists a $C^\infty$ proper Morse function $M\to\mathbf{R}$ with a unique critical point of index $0$. 

We proceed by constructing a proper, locally Lipschitz, and strictly convex function.
\begin{lem}\label{GML}
	Let $M$ be as in Theorem \ref{GM} and $p\in M$ be a fixed point. Then there exists a $C^2$ function $\chi$ such that $\chi\circ b_p: M\to \mathbf{R}$
	is a proper, locally Lipschitz, and strictly convex function.
\end{lem}

\begin{nproof}
Since $b_p$ is proper, it is bounded from below. Let $\displaystyle a = \inf_{x\in M} b_p(x)$. 
For $r \ge a$, we define 
\[
	K_0(r) := \inf \{K^\alpha(y,v)\mid y\in T_xM, b_p(x) \le r+1\}.
\]
Let 
\[
	Q(r) := \max\left\{8\left(1+\frac1{\alpha}\right)\frac1{K_0(r)}, 1\right\} \\
\]
and
\[
	\chi(t) := \int_a^t \exp\left(\int_a^s J(x)dx\right)ds
\]
where $J$ is a function $[a,\infty)\to(0,\infty)$ to be specified later.
So $\chi$ is $C^2$ on $[a,\infty)$ and 
\begin{enumerate}
	\item[(a)] $\chi'(r) \ge 1$ for $r\ge a$;
	\item[(b)] $\chi''(r) = J(r)\chi'(r)$ for $r \ge a$.
\end{enumerate}

The condition (a) above, together with the fact that $b_p$ is proper and locally Lipschitz, shows that $\chi\circ b_p$ is proper and locally Lipschitz. 

To show that $\chi\circ b_p$ is strictly convex, fix a point $q\in M$ with  $b_p(q)=r $. 
As in lemma \ref{georay} take a geodesic ray  $c_q(t): [0,\infty)\to M$  from $q$ such  that for any $t > 0$,
\[ b^{q, t}_p (x) = b_p(q)+t - d(x, c_q(t)) \]
which supports $b_p(x)$ at $q$ and $b_p(c_q(t))= t+ b_p(q)$. 
Let $v\in T_pM$ be a unit tangent vector, and $\gamma:(-\varepsilon, \varepsilon)\to M$ be a normal geodesic with $\gamma'(0) = v$. 

We construct a variation of $c_q|_{[0,Q(r)]}$
as in lemma \ref{SECVAR} with $c = c_q$, $l_0 = 1$, $l = Q(r)$ and $j(t) = 1-\frac{t}{Q(r)}$, define a function $f$ along the geodesic $\gamma$ by
\[    
    f\circ\gamma(s) := b_p(q) + Q(r) - L(s).
\]
Since
\[   
	L(s) \geq d(\gamma(s), c_q(Q(r))
\]
we have
\[
	f(x) \leq  b_p(q)+Q(r) - d(x, c_q(Q(r))) = b^{q, Q(r)} ( x)  \leq  b_p(x)
\]
for all $x$ on the geodesic $\gamma$. 
With a little abuse of notation we will denote 
\[
	df(v) = \left.\frac{\partial}{\partial s}f(\gamma(s))\right|_{s=0}, \quad H^2f(v) = \left.\frac{\partial^2}{\partial s^2}f(\gamma(s))\right|_{s=0}
\]
and similarly for $\chi\circ f$.

Put $\tau = \abs{g_{c'_q(0)}(\gamma'(0), c'_q(0))}$ so $\abs{df(v)} = \abs{L'(0)} = \tau$ and 
\[\begin{aligned}
	H^2f(v) = -L''(0) &\ge -\frac{\beta^2}{Q(r)}\left(1+\frac1\alpha\right) + \beta^2\frac{K_0(r)}{4} - \left[T_{c_q'(0)}(\gamma'(0)) - T_{c_q'(0)}(\gamma'(0)^\bot)\right]\\
	&\ge \beta^2\frac{K_0(r)}{8} - \left[T_{c_q'(0)}(\gamma'(0)) - T_{c_q'(0)}(\gamma'(0)^\bot)\right]
\end{aligned}\]
by our choice of $Q$. Choose a continuous increasing function $\lambda:[a,\infty)\to[1,\infty)$ so that 
\[
	\frac1{\lambda(r)} \sqrt{g_y(v,v)} \le F(v) \le \lambda(r)\sqrt{g_y(v,v)}
\]
for all $x\in b_p^{-1}[a,r]$, $y\in T_xM\setminus\{0\}$ and $v\in T_xM$. We have 
\[
	\frac1{\lambda(r)^2} \le \beta^2 + \tau^2 \le \lambda(r)^2,
\]
and 
\[
	H^2f(v) \ge \left(\frac1{\lambda(r)^2} - \tau^2\right)\frac{K_0(r)}{8} - \left[T_{c_q'(0)}(\gamma'(0)) - T_{c_q'(0)}(\gamma'(0)^\bot)\right].
\]

As in the proof of lemma \ref{TOTCONV} we choose a continuous functions $\tilde{\tau}:[a,\infty)\to(0,\infty)$ such that for all $q\in b_p^{-1}[a,r]$ and $\tau \le \tilde{\tau}(r)$, we have
$\frac{K_0(r)\tau^2}{8} < \frac{K_0(r)}{32\lambda(r)^2}$ and 
\[
	\abs{T_{c_q'(0)}(\gamma'(0)) - T_{c_q'(0)}(\gamma'(0)^\bot)} < \frac{K_0(r)}{32\lambda(r)^2}.
\]
Then 
\[
	H^2f(v) \ge  \frac{K_0(r)}{16\lambda(r)^2}
\]
when $\tau \le \tilde{\tau}(r)$.

Finally, let $T_0(r) := \displaystyle\sup_{x\in b_p^{-1}[a,r], y,v\in T_xM, F(y) = 1, F(v) \le \lambda(r)^2} \abs{T_y(v)}$, then $T_0$ is continuous on $[a,\infty)$, and 
\[
	H^2f(v) \ge \beta^2\frac{K_0(r)}{8} - 2T_0(r). 
\]
Put $J(r) := \frac{1}{\tilde{\tau}(r)^2}\left(2T_0(r) + \frac{K_0(r)}{16\lambda(r)^2}\right)$, then we have 
\[\begin{aligned}
	H^2(\chi\circ f)(v) &= \chi'(r)\left[H^2f(v) + \frac{df(v)^2}{\tilde{\tau}(r)^2}\left(2T_0(r) + \frac{K_0(r)}{16\lambda(r)^2}\right)\right] 
	\\&= \chi'(r) \left[H^2f(v) +\frac{\tau^2}{\tilde{\tau}(r)^2}\left(2T_0(r) + \frac{K_0(r)}{16\lambda(r)^2}\right)\right] 
	\ge \frac{K_0(r)}{16\lambda(r)^2}> 0,
\end{aligned}\]
hence $\chi\circ b_p$ is strictly convex on $M$ by lemma \ref{SUPPORT}.

\end{nproof}

\begin{rem} \rm
	We defined the weighted flag curvature with the coefficent of the $\dot{T}$ term being $1$, which turned out to be crucial in the proof. 
	This is a consequence of the fact that the combination $K(y,u) + \dot{T}_y(u)$ plays an important role in the second variation when $u$ is orthogonal to $y$ with respect to $g_y$.
\end{rem}

\bigskip
\begin{nproof}[Proof of Theorem \ref{GM}.]
    We have from the previous lemma that $\chi\circ b_p$ is proper, locally Lipschitz, and strictly convex on $M$.
	It follows from Theorem \ref{Smooth} that there is a proper, strictly convex Morse function $Z$ on $M$. 
    Since $Z$ is strictly convex along any geodesic, it has a unique critical point of index 0. Thus Theorem \ref{GM}
	follows from the standard Morse theory. 
\end{nproof}
\bigskip

\section{Proof of Theorem \ref{GM2}}

In this section we prove Theorem \ref{GM2}, a partial extension of Wu's\cite{WuH} and Sha's\cite{Sha} result about $p$-convex Riemannian manifolds.
\bigskip

\begin{nproof}[Proof of Theorem \ref{GM2}.]
	Choose a small $\delta > 0$ so that $B^-(\partial M,\delta) := \{x\in M\mid d(x,\partial M) < \delta\}$ is a collar 
	neighborhood of $\partial M$ in $M$. It follows that $M' := M\setminus B^-\left(\partial M,\frac{\delta}2\right) \cong M$ 
	and it now suffices to show that the interior of $M'$ is diffeomorphic to $\mathbf{R}^n$ hence to a Euclidean ball. 
	We shall now construct a proper, strictly convex function on the interior of $M'$ from the function $f(x) := -d(x,\partial M)$, with the same 
	idea as in the proof of Theorem \ref{GM}.
	
	Let $x\in M'$ be an interior point, $v\in T_xM$ be a unit tangent vector, and $\gamma: (-\varepsilon, \varepsilon) \to M'$ be a normal geodesic with $\gamma'(0) = v$. 
	Take $c:[0,l]\to M$ to be a minimal normal geodesic connecting $x$ to $\partial M$.
	Let $\tilde\tau > 0$ be a constant which we will specify later. 
	
	In the case when $\tau := \abs{g_{c'(0)}(\gamma'(0),c'(0))} \le \tilde\tau$, 
	we construct a variation of $c$ as in lemma \ref{SECVAR} with $l_0 = 0$, $j(t)\equiv 1$ and $c_0$ being a geodesic of $\partial M$ with $c_0(0) = c(l)$ and $c_0'(0) = E(l)$. 
	Then the function $s\mapsto -L(s)$ supports $f\circ \gamma$ at $0$, and 
	\[
		-L''(0) \ge - F(V(l))^2g_{c'(l)}(\kappa_l(0),c'(l)) - \left[T_{c'(0)}{\gamma'(0)} - T_{c'(0)}(\gamma'(0)^\bot)\right].
	\]
	Observe that $c'(l)$ is exactly the outward-pointing normal vector $\mathbf{n}$ of $\partial M$ at $c(l)$, 
	so the normal curvature of $\partial M$ is given by $\Lambda_\mathbf{n}(c_0'(0)) = - F(c_0'(0))^2 g_{c'(l)}(\kappa_l(0), c'(l))$. 
	Since $\partial M$ is strictly convex and $g_{c'(l)}(c_0'(0),c_0'(0)) = g_{c'(0)}(\gamma'(0)^\bot,\gamma'(0)^\bot)$, we may assume that $-F(c_0'(0))^2g_{c'(l)}(\kappa_l(0),c'(l)) \ge \varepsilon_0 > 0$ for some fixed 
	constant $\varepsilon_0$ depending on the manifold $M$ and $\tilde\tau$ only, hence 
	\[
		-L''(0) \ge \varepsilon_0 - \left[T_{c'(0)}{\gamma'(0)} - T_{c'(0)}(\gamma'(0)^\bot)\right].
	\]
	Since $M$ is compact, by possibly choosing a smaller $\tilde\tau$ we may assume that given any $x\in M$ and $y,v\in T_xM$ being unit tangent 
	vectors satisfying $g_y(y,v) \le \tilde\tau$, we have $\abs{T_y(v) - T_y(v^\bot)} < \frac{\varepsilon_0}{2}$, 
	where $v^\bot = v - g_y(y,v)y$. It follows that 
	\[
		-L''(0) \ge \frac{\varepsilon_0}{2} > 0.
	\]
	
	On the other hand, if $\tau > \tilde\tau$, 
	we construct a variation of $c$ as in lemma \ref{SECVAR} with $l_0 = 0$, $j(t) = 1 - \frac{t}{l}$ and $c_0(s)\equiv c(l)$. 
	The function $s\mapsto -L(s)$ supports $f\circ \gamma$ at $0$, and 
	\[
		-L'(0) = g_{c'(0)}(\gamma'(0),c'(0)),
	\]
	\[
		-L''(0) \ge - \left(1+\frac1\alpha\right)\frac{\beta^2}{l} - \left[T_{c'(0)}{\gamma'(0)} - T_{c'(0)}(\gamma'(0)^\bot)\right].
	\]

	Choose $\lambda \ge 1$ so that 
	\[
		\frac1{\lambda}\sqrt{g_y(v,v)} \le F(v) \le \lambda \sqrt{g_y(v,v)}
	\]
	for any $x\in M$, $y\in T_xM\setminus\{0\}$ and $v\in T_xM$. Now we have 
	\[
		\frac1{\lambda^2} \le \beta^2 + \tau^2 \le \lambda^2,
	\]
	and 
	\[
		F(\gamma'(0)^\bot) \le \beta\lambda \le \lambda^2
	\]

	By construction $l \ge \frac{\delta}{2}$. Let $T_0 = \displaystyle\sup_{x\in M, y,v\in T_xM, F(y) = 1, F(v)\le \lambda^2}\abs{T_y(v)}$, we have 
	\[
		-L''(0) \ge - \left(1+\frac1\alpha\right)\frac{2\beta^2}{\delta} - 2T_0 \ge - \left(1+\frac1\alpha\right)\frac{2\lambda^2}{\delta} - 2T_0=: -U_0.
	\]
	Now we define $\chi: (-\infty,0]\to\mathbf{R}$ by 
	\[
		\chi(x) := \int_0^x \exp\left(\frac{U_0 + \varepsilon_0}{\tilde\tau^2}t\right)dt, 
	\]
	so $\chi'(x) = \exp\left(\frac{U_0 + \varepsilon_0}{\tilde\tau^2}x\right) > 0$ and $\chi''(x) = \frac{U_0 + \varepsilon_0}{\tilde\tau^2}\chi'(x)$. Since $\chi$ is (strictly) increasing, 
	$\chi\circ (-L)$ supports 
	$\chi\circ f\circ \gamma$ at $0$ in either of the above cases. Moreover, 
	\[\begin{aligned}
		(\chi\circ (-L))''(0) =& \chi'(-L(0))\left[-L''(0) + \frac{U_0 + \varepsilon_0}{\tilde\tau^2}L'(0)^2\right] \\
		=& \chi'(-L(0))\left[-L''(0) + \frac{U_0 + \varepsilon_0}{\tilde\tau^2}\tau^2\right] \ge \chi'(-\sup_{x\in M}d(x,\partial M))\frac{\varepsilon_0}{2} > 0.
	\end{aligned}\]
	By lemma \ref{SUPPORT} $\chi\circ f$ is a strictly convex function defined on the interior of $M'$. Finally, since $f$ is proper, Lipschitz and bounded below, we have that $\chi\circ f$ is proper and Lipschitz. 
	With the same approximating argument we conclude that the interior of $M'$ is diffeomorphic to $\mathbf{R}^n$ hence to a Euclidean ball, and we are done.
\end{nproof}

\section*{Appendix}

We now sketch a proof of the smoothing theorem \ref{Smooth}.

As in \cite{GrWu2}, for each compact subset $K\subset M$, we may define a metric $d_K$ on the space of $C^\infty$ functions in a neighborhood of $K$, which is independent of the choice of
the metric on the manifold and gives the $C^\infty$ topology on the function space. 

By corollary 1 of Theorem 4.1 in \cite{GrWu2}, the proof of theorem \ref{Smooth} reduces to the following 3 lemmas:
\begin{lemmaa}\label{maxClo}
	The set of locally Lipschitz strictly convex functions has the maximum closure property, in the sense that given two locally Lipschitz strictly convex functions $f_1, f_2$, 
	then $\max(f_1,f_2)$ is	locally Lipschitz and strictly convex.
\end{lemmaa}
\begin{lemmab}\label{smoSta}
	The set of locally Lipschitz strictly convex functions has the $C^\infty$-stability property, in the sense that given any compact set $K\subset M$ and 
	locally Lipschitz strictly convex function $f$,
	then there is a positive number $\varepsilon > 0$ such that for any $C^\infty$ function $\phi$ with $d_K$ norm less than $\varepsilon$, $f+\phi$ is locally Lipschitz and strictly
	convex on a neighborhood of $K$.
\end{lemmab}
\bigskip
The above two lemmas are elementary. The last one is less trivial:
\begin{lemmac}\label{locApprox}
	The set of locally Lipschitz strictly convex functions has the local approximation property, 
    in the sense that given any $x\in M$, there is an open neighborhood $U$ of $x$ with the
	following property: 
    Let $L\subset K$ be compact sets and $V$ be an open set such that $K\subset V\subset U$. 
    Given any locally Lipschitz strictly convex function $f$ on $V$, $C^\infty$ on $L$, there exists an open neighborhood $W$ of $K$ in $V$, 
    such that for any positive constant $\varepsilon > 0$, there is a $C^\infty$, locally Lipschitz, and strictly convex function $g$ on $V$, satisfying $\sup_K\abs{g - f} < \varepsilon$ and $d_L(f,g)< \varepsilon$.
\end{lemmac}

\begin{nproof}
	We choose a Riemannian metric $g$ on $M$, and let $U$ be a pre-compact neighborhood of $x$ in $M$. 
    Now choose $W$ and $\delta$ so that the $\delta$-neighborhood of $W$, with respect to $g$, is contained in $V$.
    We further assume that $\exp^g$ maps the $\delta$-ball in $T_pM$ diffeomorphically onto its image for all $p\in W$ 
    where $\exp^g$ is the exponential map of $g$.
	For any $p\in W$, let 
	\[
		f_\delta(p) = \dfrac{1}{\delta^n}\int_{T_pM} f(\exp^g_p(v))\phi\left(\frac{\norm{v}_g}{\delta}\right)d\mu_p
	\]
	where $\phi$ is a nonnegative $C^\infty$ function supported in $[-1,1]$, constant in a neighborhood of $0$, and satisfies $\int_{\mathbf{R}^n}\phi(\norm{v})dv = 1$,
	and $d\mu_p$ is the Lebesgue measure on $T_pM$ relative to the Riemannian inner product $g$. A standard argument in Riemannian geometry shows that for sufficiently small
	$\delta$, $f_\delta$ is a well-defined $C^\infty$ function on $W$, and that $f_\delta$ converges to $f$ as $\delta \to 0$, in the $C^0$ topology on $K$ and in the $C^\infty$ topology on $L$.
	
	Fix $p\in W$, and an $F$-unit vector $v\in T_pM$, and let $\gamma: (-\varepsilon, \varepsilon)\to V$ be the $F$-geodesic with $\gamma'(0) = v$. Let $\mathrm{P}_{s}u$ be the vector in $T_{\gamma(s)}M$ obtained from $u\in T_pM$ by 
	a $g$-parallel transport along $\gamma$. Then 
	\[
		f_\delta(\gamma(-t)) + f_\delta(\gamma(t)) = \dfrac{1}{\delta^n}\int_{T_pM}\left[f(\exp^g_{\gamma(-t)}\mathrm{P}_{-t}u) + f(\exp^g_{\gamma(t)}\mathrm{P}_{t}u)\right]\phi\left(\frac{\norm{u}_g}{\delta}\right)d\mu_p
	\]
	Let $u\in T_pM$ be such that $\norm{u}_g \le \delta$, and $c_0: (-\varepsilon, \varepsilon)\to M$ be defined by $c_0(t) = \exp^g_{\gamma(t)}\mathrm{P}_tu$. Choosing $\delta$ small enough, there is a unique (not necessarily normal)
	$F$-geodesic $\gamma_u:(-\varepsilon, \varepsilon)\to M$ with $\gamma_u(0) = c_0(0)$ and $\gamma_u'(0) = c_0'(0)$. 
	
	By the $C^\infty$ dependence of the solutions of ordinary differential equation on the initial conditions, we check that 
    both $c_0$ and $\gamma_u$ converge to $\gamma$ in the $C^\infty$ topology. Then by Lemma 3 in \S{3} of \cite{GrWu1} it can be shown that for any given $\zeta > 0$,
	\[
		d^g(c_0(t), \gamma_u(t)) \le \zeta t^2
	\]
	holds for all sufficiently small $t$ and $\norm{u}_g\le \delta$. Now let $L_p$ be a $g$-Lipschitz constant of $f$ on $\overline{V}$, this implies 
	\[\begin{aligned}
		f_\delta(\gamma(-t)) + f_\delta(\gamma(t)) \ge \dfrac{1}{\delta^n}\int_{T_pM}\left[f(\gamma_u(-t)) + f(\gamma_u(t)) - 2L_p\zeta t^2\right]\phi\left(\frac{\norm{u}_g}{\delta}\right)d\mu_p \\
		= \dfrac{1}{\delta^n}\int_{T_pM}\left[f(\gamma_u(-t)) + f(\gamma_u(t)) \right]\phi\left(\frac{\norm{u}_g}{\delta}\right)d\mu_p - 2L_p\zeta t^2
	\end{aligned}\]
	Hence 
	\[\begin{aligned}
		f_\delta(\gamma(-t))& + f_\delta(\gamma(t)) - 2f_\delta(x) \\
		&\ge \dfrac{1}{\delta^n}\int_{T_pM}\left[f(\gamma_u(-t)) + f(\gamma_u(t)) - 2f(\gamma_u(0)) \right]\phi\left(\frac{\norm{u}_g}{\delta}\right)d\mu_p - 2L_p\zeta t^2
	\end{aligned}\]
	Since $f$ is strictly convex, 
	\[
		f(\gamma_u(-t)) + f(\gamma_u(t)) - 2f(\gamma_u(0)) \ge \frac{\delta_{\overline{W}}}{2}t^2 F(c'_0(0))^2
	\]
	for some $\delta_{\overline{W}} > 0$ and small enough $t$.
	Since $c'_0(0)$ depends smoothly on $u$, by taking a uniform upper bound of $F(c'_0(0))$ for $\norm{u}_g\le \delta$ on $\overline{W}$, we have that 
    \[
        \dfrac{1}{\delta^n}\int_{T_pM}\left[f(\gamma_u(-t)) + f(\gamma_u(t)) - 2f(\gamma_u(0)) \right]\phi\left(\frac{\norm{u}_g}{\delta}\right)d\mu_p \ge M_0t^2
    \]
	for some other constant $M_0 > 0$ and small enough $t$.
	Taking $\zeta = \frac{M_0}{4L_p}$ and we have that $f_\delta$ is also strictly convex for sufficiently small $\delta$.
\end{nproof}

\bibliographystyle{plain}

\vspace{0.6cm}

\noindent Zhongmin Shen

\noindent Department of Mathematical Sciences, Indiana
University-Purdue University Indianapolis, IN 46202-3216, USA.

\noindent \verb"zshen@math.iupui.edu"
\bigskip

\noindent Runzhong Zhao

\noindent School of Mathematical Sciences, Xiamen University, 361005, China.

\noindent \verb"runzzhao@xmu.edu.cn"

\end{document}